\documentstyle[12pt]{amsart}

\makeatletter %restore definition changed by amsart
\def\cases#1{\left\{\,\vcenter{\normalbaselines\m@th
    \ialign{$##\hfil$&\quad##\hfil\crcr#1\crcr}}\right.}
\makeatother
\let\<\langle
\let\>\rangle

\author{Jie Wu}
\address{MSRI, 1000 Centennial Drive, Berkeley, CA 94720, USA}
\email{jwu@@msri.org}
\title{A Product Decomposition of $\Omega^3_0\Sigma\bold RP^2$}
\thanks{Research at MSRI is supported in part by NSF grant DMS-9022140.}

\newtheorem{theorem}{Theorem}[section]

\newtheorem{lemma}[theorem]{Lemma}

\newtheorem{corollary}[theorem]{Corollary}

\begin{document}
\maketitle
\begin{abstract}
We give a specific product decomposition of the base-point path connected component of the triple loop
space of the suspension of the projective plane.
\end{abstract}
\section{Introduction}
   In this paper, we give a specific product decompositions of the base-point 
 path connected component of the triple loop space of the suspension of the 
 projective plane, which is denoted by $\Omega^3_0\Sigma\bold RP^2$. Let
 $\bold RP_a^b=\bold RP^b/\bold RP^{a-1}$, let $X<n>$ be the
 $n$-connected covering of the space $X$ and let $P^n(2)=\Sigma^{n-2}\bold RP^2$,
the $n$-dimensional mod $2$ Moore space. In particular, $P^3(2)=\Sigma\bold RP^2$. Our theorem is as follows.
  
\begin{theorem}
There is a homotopy equivalence
$$
\Omega^3_0P^3(2)\simeq\Omega^2(S^3\<3\>)\times\Omega^3_0(P^6(2)\vee\Sigma\bold RP^4_2)
$$ 
localized at $2$.
\end{theorem}

This theorem was first given in my thesis where we only gave a product decomposition of
$\Omega^4_0(\Sigma\bold RP^2)$. The author would like to thank his thesis advisor Fred Cohen for pointing
out that our original decomposition for $\Omega^4_0(\Sigma\bold RP^2)$ can be delooping. He showed that
the power map 
$$
4:\Omega^2(S^3\<3\>)\rightarrow\Omega^2(S^3\<3\>)
$$
is null localized at $2$, his proof will be used
in this paper. The application of this theorem to calculate
the homotopy groups $\pi_r(P^3(2))$ for $r\leq 11$ was given in my thesis
[Wu]. An independent calculation for $\pi_*(P^3(2))$ can be found out in [Mu]. The homotopy groups of
$\Sigma\bold RP^2$ may be useful in low dimensional topology.

\section{Proof of the Theorem 1.1}

In this section, all spaces are localized at $2$ and $H_*(X)$ means the homology with coefficients in
field $\bold F_2$, the field with two elements.  Let
$j:P^3(2)\rightarrow BSO(3)$ be the inclusion of the $3$-skeleton of $BSO(3)$ and let $X_2$ be the pull
back in the diagram

$$
\hspace{1.5in}
\begin{array}{cccccc}
& BSO(2) &\rightarrow & BSO(3) \\  
& \uparrow &  & \uparrow  \\  
& X_2 &\rightarrow &P^3(2), \\  
\end{array}
$$
where the map $BSO(2)\rightarrow BSO(3)$ is induced by the canonical injection $SO(2)\rightarrow SO(3)$.
Thus there is a fibration 
$$
S^2\rightarrow X_2 \rightarrow P^3(2).
$$
Let $i:S^2\rightarrow X_2$ be the inclusion of the fibre.

\begin{lemma}
The map $\Omega^2(i):\Omega^2(S^2\<3\>)\rightarrow\Omega^2(X_2\<3\>)$ is null.
\end{lemma}
\noindent{\em Proof:} Consider the commutative diagram

$$
\hspace{0.5in}
\begin{array}{cccccccccccc}
&\rightarrow&\pi_3(BSO(3))          &\stackrel{0}{\rightarrow}       &\pi_2(S^2)  &\hookrightarrow               &\pi_2(BSO(2))&\rightarrow&\pi_2(BSO(3))           &\rightarrow0\\  
&           & \uparrow              &                                & \Vert      &                              &\uparrow     &           & \uparrow               &            \\  
&\rightarrow&\pi_3(P^3(2))&\stackrel{0}{\rightarrow}      &\pi_2(S^2)   &\stackrel{i_*}{\hookrightarrow}&\pi_2(X_2)  &\rightarrow
&\pi_2(P^3(2))&\rightarrow0\\  
\end{array}
$$
where the top and bottom rows are long exact sequence of homotopy groups. Thus
$\pi_1(X_2)=0$, $\pi_2(X_2)\cong\pi_2(BSO(2))=\pi_2(\bold CP^{\infty})\cong\bold Z$
and $i_*:\pi_2(S^2)\rightarrow\pi_2(X_2)$ is of degree $2$. Let $\phi:S^2\rightarrow X_2$ be a 
representive of the generator in $\pi_2(X_2)$. Then there is a homotopy commutative diagram
$$
\hspace{1.5in}
\begin{array}{cccccc}
& S^2 & \stackrel{i}{\rightarrow} & X_2 \\  
& \Vert &  & \uparrow \phi \\  
& S^2&\stackrel{[2]}{\rightarrow} & S^2, \\  
\end{array}
$$
where $[2]:S^2\rightarrow S^2$ is of degree $2$. Now consider the homotopy commutative diagram
$$
\hspace{1.5in}
\begin{array}{cccccc}
& \Omega S^2 & \stackrel{H_2}{\rightarrow} & \Omega S^3 \\  
&\Omega([2])\downarrow &  & \downarrow \Omega([4]) \\  
&\Omega S^2&\stackrel{H_2}{\rightarrow} & \Omega S^3, \\  
\end{array}
$$
where $H_2$ is the second James-Hopf map. Notice that $S^3$ is an H-space. We have
$$
\Omega([4])\simeq 4:\Omega S^3\rightarrow \Omega S^3,
$$
where $4:\Omega S^3\rightarrow \Omega S^3$ is the power map of degree $4$. 
Consider the EHP sequence $S^3\rightarrow\Omega(S^3\<3\>)\stackrel{H_2}{\rightarrow}\Omega S^5$, the map
$S^3\rightarrow\Omega(S^3\<3\>)$ is a representative of the nontrivial element $\eta$ in
 $\pi_4(S^3)\cong\bold Z/2$. Notice that the composite 
 $$
 S^3\stackrel{\eta}{\rightarrow}\Omega(S^3\<3\>)\stackrel{2}{\rightarrow}\Omega(S^3\<3\>)
 $$
 is null. Thus the map $2\simeq\Omega(2):\Omega^2(S^3\<3\>)\rightarrow\Omega^2(S^3\<3\>)$ factors through
$\Omega(H_2):\Omega^2(S^3\<3\>)\rightarrow\Omega^2S^5$ up to homotopy. By a result of Cohen [C1, or C2], the map
$$
2\cdot\Omega(H_2):\Omega^2(S^3\<3\>)\rightarrow\Omega^2S^5
$$
is null. Thus the power map 
$$
4:\Omega^2(S^3\<3\>)\rightarrow\Omega^2(S^3\<3\>)
$$
is null.
Notice that the map
$\Omega (H_2):\Omega^2(S^2\<3\>)\rightarrow \Omega^2(S^3\<3\>)$ is a homotopy equivalence. Thus
$\Omega^2([2]):\Omega^2(S^2\<3\>)\rightarrow \Omega^2(S^2\<3\>)$ is null and so the boundary map
$\Omega^2(i):\Omega^2(S^2\<3\>)\rightarrow \Omega^2(X_2\<3\>)$ is null, which is the assertion.

\begin{corollary}
There is a homotopy equivalence
$$
\Omega^3_0P^3(2)\simeq\Omega^2(S^3\<3\>)\times\Omega^3_0(X_2).
$$
\end{corollary}

\begin{lemma}
The (mod $2$) homology of $X_2\<2\>$ is as follows:
$$
{\bar H}_r(X_2\<2\>)=\cases{\bold Z/2&r=3,4,6\cr \bold Z/2\oplus\bold Z/2&r=5\cr 0&otherwise.}
$$
\end{lemma}
\noindent{\em Proof:}
Consider the homotopy commutative of fibre sequences
$$
\hspace{0.5in}
\begin{array}{cccccccccc}
& S^2 & \rightarrow&BSO(2)=\bold CP^{\infty}&\rightarrow &BSO(3) \\  
&\Vert&            &\uparrow               &             & \uparrow \\  
&S^2  &\rightarrow &X_2                    &\rightarrow   & P^3(2) \\
&\uparrow&         &\uparrow               &              &\uparrow\\
&*    &\rightarrow &X_2\<2\>                 &\stackrel{=}{\rightarrow}&X_2\<2\>\\ 
\end{array}
$$
Thus the homotopy fibre of $X_2\<2\>\rightarrow P^3(2)$ is $SO(3)$.
The assertion follows by inspecting the Serre spectral sequence for the principle fibre sequence
$$
SO(3)\rightarrow X_2\<2\>\rightarrow P^3(2).
$$ 
\par
\par
Recall that$ H_*(\Omega P^3(2))\cong T(u,v)$
as algebras with $dim(v)=2$ and $Sq^1_*v=u$, where $T(u,v)$ is the tensor algebra generated by $u$ and
$v$.

\begin{lemma}
$H_*(\Omega (X_2\<2\>))\cong T(V)$ as subalgebras of $H_*(\Omega P^3(2))$, where $V$ has a basis
$$
\{u^2,[u,v],v^2+u[u,v],[u^2,v],[v^2,u]+u[u^2,v]\}
$$
in $T(u,v)\cong H_*(\Omega P^3(2))$ and $[x,y]$ is the commutator of $x$ and $y$.
\end{lemma}
\noindent{\em Proof:} Notice that $\Omega j_*:H_*(\Omega P^3(2))\rightarrow H_*(SO(3))$ is an epimorphism and
$$\pi_1(\Omega P^3(2))\cong\pi_1(SO(3)).$$ Thus the Serre spectral sequence for the fibre sequence
$$
\hspace{0.5in}
\begin{array}{cccccc}
&\Omega (X_2\<2\>) & \stackrel{\Omega i}{\rightarrow} &\Omega P^3(2)&\stackrel{\Omega j}{\rightarrow} & SO(3) \\  
\end{array}
$$
collapses and therefore there is a short exact sequence of Hopf algebras
$$
1\rightarrow H_*(\Omega (X_2\<2\>))\rightarrow H_*(\Omega P^3(2))\rightarrow H_*(SO(3))\rightarrow1.
$$
Thus 
$H_*(\Omega (X_2\<2\>))$ is the cotensor product of $\bold F_2$ and $H_*(\Omega P^3(2))$ over
$H_*(SO(3))$ ( see [MM] ).
\par
 We first show that $V\subseteq H_*(\Omega (X_2\<2\>))$. Let $\psi$ be the comultiplication
in the Hopf algebra $H_*(\Omega P^3(2))\cong T(u,v)$. Then $\psi$ is determined by the generators.
$$
\psi(u)=u\otimes1+1\otimes u
$$
and
$$
\psi(v)=v\otimes1+1\otimes v+u\otimes u.
$$
The comultiplication $\psi$ on the basis for $V$ is as follows:
$$
\psi(u^2)=u^2\otimes 1+1\otimes u^2,
$$
$$
\psi([u,v])=[u,v]\otimes1+1\otimes[u,v],
$$
$$
\psi(v^2+u[u,v])=(v^2+u[u,v])\otimes1+1\otimes(v^2+u[u,v])+u^2\otimes u^2,
$$
$$
\psi([u^2,v])=[u^2,v]\otimes1+1\otimes[u^2,v],
$$
and
$$
\psi([u,v^2]+u[u^2,v])=([u,v^2]+u[u^2,v])\otimes1+1\otimes([u,v^2]+u[u^2,v]).
$$
\par
Notice that $H_*(SO(3))\cong E(u,v)$, the exterior algebra generated by $u$ and $v$. Thus
$V\subseteq H_*(\Omega (X_2\<2\>))$. Let $\phi:T(V)\rightarrow H_*(\Omega (X_2\<2\>))$ be the algebra map
induced by the inclusion $V\subseteq H_*(\Omega (X_2\<2\>))$ and let $A=\Omega i_*\circ\phi(T(V))$. Then
$A$ is the subalgebra of $H_*(\Omega P^3(2))$ generated by $V$, which is also a tensor algebra.
Notice that $V\rightarrow Q(A)$ is an isomorphism. Thus $A\cong T(V)$ ( see [CMN, Proposition 3.9]) and
so
$\phi:T(V)\rightarrow H_*(\Omega (X_2\<2\>))$ is a monomorphism. Notice that
Poincare series
$$
\chi(H_*(\Omega (X_2\<2\>)))=\chi(H_*(\Omega P^3(2)))/\chi(H_*(SO(3)))
$$
$$
=(1-t-t^2)^{-1}(1+t+t^2+t^3)^{-1}=(1-t^2-t^3-2t^4-t^5)^{-1}=\chi(T(V)).
$$
Thus $\phi:T(V)\rightarrow H_*(\Omega (X_2\<2\>))$ is an isomorphism, which is the assertion.

\begin{lemma}
The homology suspension $\sigma_*:Q{\bar H}_*(\Omega (X_2\<2\>))\rightarrow {\bar H}_{*+1}(X_2\<2\>)$ is an isomorphism.
\end{lemma}
\noindent{\em Proof:} Let $x_2,x_3,x_4,y_4$ and $y_5$ be in $H_*(\Omega(X_2\<2\>))$ corresponding to
$u^2,[u,v],v^2+u[u,v],[u^2,v]$ and $[v^2,u]+u[u^2,v]$ in $T(u,v)$, respectively. Consider mod $2$ Serre spectral
sequence $E^r_{*,*}$ for the principle fibre sequence 
$$
\Omega (X_2\<2\>)\rightarrow *\rightarrow X_2\<2\>.
$$
Notice that, 
for $*\leq5$, ${\bar H}_*(\Omega (X_2\<2\>))$ has a basis
$$
x_2,x_3,x_2^2,y_4,x_4,[x_2,x_3],x_2x_3,y_5.
$$
Let $a_{j+1}=\sigma_*(x_j)$ and $b_{j+1}=\sigma_*(y_j)$. 
Then 
$$
d_3(a_3)=x_2,
$$
$$
d_3(a_3x_2)=x_2^2,
$$
$$
d_3(a_3x_3)=x_2x_3,
$$
$$
d_4(a_4)=x_3,
$$
$$
d_4(a_4x_2)=x_3x_2,
$$
and
$$
E^r_{s,t}=0
$$
for $r\geq 6$,$s>0$,$t>0$ and $s+t\leq 6$. Thus there exists $a_5,b_5,b_6\in E^2_{*,0}$ such that
$d_5(a_5)=x_4$,$d_5(b_5)=y_4$ and $d_6(b_6)=y_5$. The assertion follows.

\begin{lemma} The element
$$x_3=[u,v]$$
in $H_3(\Omega (X_2\<2\>))$ is spherical.
\end{lemma}
\noindent{\em Proof:} Consider the exact sequence of homotopy groups
$$
\pi_4(P^3(2))\rightarrow\pi_3(SO(3))\rightarrow\pi_3(X_2\<2\>)\rightarrow\pi_3(P^3(2))\rightarrow\pi_2(SO(3))=0.
$$
Notice that $\pi_3(SO(3))=\bold Z$ and $\pi_4(P^3(2))$ is torsion. Thus $\pi_3(X_2\<2\>)=\bold Z$ or $\bold Z\oplus G$
for some $2$-torsion abelian group $G$. Now $X_2\<2\>$ is 2-connected and 
$\bold F_2=H_3(X_2\<2\>;\bold F_2)\cong \pi_3(X_2\<2\>)\otimes\bold F_2$. Thus $\pi_3(X_2\<2\>)=\bold Z$ and
so $\pi_2(\Omega (X_2\<2\>))=\bold Z$. Now the $3$-skeleton $(\Omega (X_2\<2\>))_{(3)}$ of $\Omega (X_2\<2\>)$ is 
the homotopy cofibre of an attaching map $\phi:S^2\rightarrow S^2$. Notice that
$$
H_2((\Omega (X_2\<2\>))_{(3)};\bold Z)\cong H_2(\Omega (X_2\<2\>);\bold Z)=\bold Z.
$$
Thus $\phi$ is of degree $0$ and so the $3$-skeleton $(\Omega (X_2\<2\>))_{(3)}\simeq S^2\vee S^3$. 
The assertion follows.\\
\\
\par

\noindent{\em Proof of Theorem 1.1:}
By Corollary 2.2, it suffices to show that
$$
X_2\<2\>\simeq P^6(2)\vee\Sigma\bold RP^4_2.
$$
We use the notation in Lemmas 2.4 and 2.5. By Lemma 2.4, Steenrod operations on ${\bar H}_*(X_2\<2\>)$ are given by
$$
Sq^2_*a_5=a_3,Sq^1_*a_5=a_4,Sq^1_*b_6=b_5.
$$
Thus
$$
{\bar H}_*(X_2\<2\>)\cong M\oplus N
$$
as modules over Steenrod algebra, where $M$ is generated by $a_3,a_4,a_5$ and $N$ is generated by $b_5,b_6$.
By Lemma 2.6, the element $x_3$ is spherical in $H_3(\Omega (X_2\<2\>))$. Let $f:S^3\rightarrow \Omega
(X_2\<2\>)$ so that 
$f_*(\iota_3)=x_3$, where $\iota_3$ is a generator for $H_3(S^3)$. Consider the relative Samelson
product for the multiplicative fibre sequence
$$
\Omega (X_2\<2\>)\rightarrow\Omega P^3(2)\rightarrow SO(3)
$$
( see [N1] ). There is a relative Samelson product
$$
[f,E]:S^3\wedge\bold RP^2=P^5(2)\rightarrow\Omega (X_2\<2\>),
$$
where $E:\bold RP^2\rightarrow \Omega P^3(2)$ is the suspension.
Notice that both $\iota_3$ and $u$ are primitive. Then
$$
[f,E]_*(\iota_3\otimes u)=[f_*(\iota_3),u]=[[u,v],u]=[u^2,v]
$$
in $T(u,v)$. Let $\theta=\sigma\circ\Sigma[f,E]$ be the composite
$P^6(2)\rightarrow\Sigma\Omega (X_2\<2\>)\stackrel{\sigma}{\rightarrow} X_2\<2\>$ and let $u_5$ and $v_6$ be the bases for $H_5(P^6(2))$
and $H_6(P^6(2))$, repectively. Then $\theta_*(u_5)=\sigma_*([u^2,v])=b_5$ and
$Sq^1_*\theta_*(v_6)=b_5$. Thus $\theta_*(v_6)=b_6$ which is the only nonzero element in $H_6(X_2\<2\>)$
and so 
$$
\theta_*:H_*(P^6(2))\rightarrow H_*( X_2\<2\>)
$$
is one to one. Let $B^5$ be the homotopy cofibre of
$\theta:P^6(2)\rightarrow X_2\<2\>$. Then  ${\bar H}_*(B^5)$ has a basis
$\{a_3,a_4,a_5\}$. Notice that
$a_4=\sigma_*(x_3)$ is spherical. Thus $B^5$ is the homotopy cofibre of an attaching map
$\phi:S^4\rightarrow S^3\vee S^4=(X_2\<2\>)_{(4)}$. Notice that the homotopy fibre of the pinch map
$X_2\<2\>\rightarrow B^5$ is $4$-connected and $\pi_4(X_2\<2\>)\cong\pi_4(B^5)$. Thus the composite
$$
S^4\stackrel{\phi}{\rightarrow}S^4\vee S^3\rightarrow X_2\<2\>
$$
is null, or the injection $S^4\vee S^3\rightarrow X_2\<2\>$ extends to a map $s:B^5\rightarrow X_2\<2\>$. Now the
composite
$$
P^6(2)\vee B^5\stackrel{\theta\vee s}{\rightarrow}X_2\<2\>\vee
X_2\<2\>\stackrel{\bigtriangledown}{\rightarrow} X_2\<2\>
$$
is a homotopy equivalence, where $\bigtriangledown$ is the fold map.  

Consider homotopy commutative diagram of cofibre sequences

$$
\hspace{0.5in}
\begin{array}{cccccccccc}
& S^4 & \stackrel{\phi}{\rightarrow}&S^3\vee S^4     &\rightarrow& B^5        & \rightarrow &S^5 \\  
&\Vert&                             &p\downarrow &           & \downarrow &             &\Vert \\  
&S^4 &\stackrel{[k]}{\rightarrow}   & S^4            &\rightarrow&P^5(k)      &\rightarrow  &S^5, \\  
\end{array}
$$
where $p$ is the projection.
Notice that $Sq^1_*a_5=a_4$. Thus $P^5(k)\simeq P^5(2)$ localized at $2$ and $k=2k'$ with $k'\equiv 1(2)$.
Notice that $Sq^2_*a_5=a_3$. Thus $B^5$ is not homotopy equivalent to $P^5(2)\vee S^3$ and therefore the composite

$$
\hspace{0.5in}
\begin{array}{cccccccccc}
&S^4 & \stackrel{\phi}{\rightarrow}&S^3\vee S^4&\stackrel{p}{\rightarrow} &S^3 \\  
\end{array}
$$
is essential which is equal to $\eta$ in $\pi_4(S^3)=\bold Z/2$,
where $p$ is the projection. Thus the homotopy class
$[\phi]=\eta+2k'[\iota]$ in $\pi_4(S^3\vee S^4)\cong\bold Z/2\oplus\bold Z$ with generators $\eta$
and $[\iota]$, where $k'\equiv 1(2)$. Notice that $\Sigma\bold RP^4_2$ is the homotopy cofibre of
$(\eta,[2]):S^4\rightarrow S^3\vee S^4$. Thus $B^5\simeq\Sigma\bold RP^4_2$ and the assertion follows.


\begin{thebibliography}{Lellmannaa}
\bibitem[C1]{C1} F. R. Cohen, {\em Two-Primary analogues of Selick's theorem and the Kahn-Priddy theorem for the $3$-sphere, } Topology, 23(1984),401-421.
\bibitem[C2]{C2} F. R. Cohen, {\em A course in some aspects of classical homotopy theory, } S.L.M 1286(1986),1-92.
\bibitem[CMN]{CMN} F. R. Cohen, J. Moore and J. Neisendorfer, {\em Torsion in homotopy groups,}  Ann. of Math. 109(1979),121-168.
\bibitem[Mu]{Mu} J. Mukai, {\em Some homotopy groups of $M^3$,} preprint.
\bibitem[MM]{MM} J. Milnor and J. Moore, {\em On the structure of Hopf Algebras,  } Ann. of Math. 81 (1965), 211-264.
\bibitem[N1]{N1} J. Neisendorfer, {\em Primary homotopy theory}, Memoirs. AMS, 25 (1980).
\bibitem[Wu]{Wu} J. Wu, Thesis, University of Rochester (1995).
\end{thebibliography}
\end{document}